\documentclass[conference]{IEEEtran}
\IEEEoverridecommandlockouts
\usepackage{cite}
\usepackage{amsmath,amssymb,amsfonts}
\usepackage{graphicx}
\usepackage{textcomp}
\usepackage{lipsum}
\usepackage{overpic}
\usepackage{booktabs}
\usepackage{xcolor}
\usepackage{math_operators}
\usepackage{subcaption}
\usepackage{siunitx}
\usepackage{url}

\makeatletter
\newcommand{\linebreakand}{%
  \end{@IEEEauthorhalign}
  \hfill\mbox{}\par
  \mbox{}\hfill\begin{@IEEEauthorhalign}
}
\makeatother

\def\BibTeX{{\rm B\kern-.05em{\sc i\kern-.025em b}\kern-.08em
    T\kern-.1667em\lower.7ex\hbox{E}\kern-.125emX}}
\begin{document}

\title{Solving Subset Sum Problems using Quantum Inspired Optimization Algorithms with Applications in Auditing and Financial Data Analysis
}

\author{
\IEEEauthorblockN{1\textsuperscript{st} David Biesner}
\IEEEauthorblockA{\textit{Fraunhofer IAIS and University of Bonn} \\
Sankt Augustin and Bonn, Germany \\
% \href{mailto:david.biesner@iais.fraunhofer.de}{david.biesner@iais.fraunhofer.de}}
\url{david.biesner@iais.fraunhofer.de}}
\and
\IEEEauthorblockN{2\textsuperscript{nd} Thore Gerlach}
\IEEEauthorblockA{\textit{Fraunhofer IAIS and University of Bonn} \\
Sankt Augustin and Bonn, Germany}
\linebreakand
\IEEEauthorblockN{4\textsuperscript{th} Christian Bauckhage}
\IEEEauthorblockA{\textit{Fraunhofer IAIS and University of Bonn} \\
Sankt Augustin and Bonn, Germany}
\and 
\IEEEauthorblockN{5\textsuperscript{th} Bernd Kliem}
\IEEEauthorblockA{\textit{PWC GmbH WPG} \\
Frankfurt a. M., Germany}
\and
\IEEEauthorblockN{3\textsuperscript{rd} Rafet Sifa}
\IEEEauthorblockA{\textit{Fraunhofer IAIS} \\
Sankt Augustin, Germany}
}
% \author{Anonymous Authors\vspace{3.5cm}}

\maketitle

\begin{abstract}
Many applications in automated auditing and the analysis and consistency check of financial documents
can be formulated in part as the subset sum problem:
Given a set of numbers and a target sum, find the subset of numbers that sums up to the target.
The problem is NP-hard and classical solving algorithms are therefore not practical to use in many real applications.

We tackle the problem as a QUBO (quadratic unconstrained binary optimization) problem and show how gradient descent on Hopfield Networks reliably finds solutions for both artificial and real data.
We outline how this algorithm can be applied by adiabatic quantum computers (quantum annealers) and specialized hardware (field programmable gate arrays) for digital annealing and run experiments on quantum annealing hardware.\footnote{To be published in proceedings of IEEE International Conference on     Machine Learning Applications IEEE ICMLA 2022.}
\end{abstract}


\section{Introduction and Problem Statement}

The financial auditing process involves the writing and proofreading of financial reports for the audited company.
This process is still largely a manual one:
auditors must read documents, compare document content with previous reports, check the completeness of the content to financial regulation checklists and check against both compliance and mathematical errors.

One aspect of mathematical correctness is the correctness of numerical tables,
for example describing profit and loss for a given year or quarter.
All values in these tables must of course correspond to the actual financial situation of the company,
which includes the correctness of sums in the tables.
For example, for a table depicting the revenue, expenses and income,
the values for expenses and income must sum up to the revenue.
Audited reports are provided most commonly as .pdf-files, with no machine readable information on the structure of the sums in the tables available. 
Even though the numbers in the tables and the corresponding calculations are likely to stem from a table processing software (e.g. Excel) and  copied directly, mistakes in the process are possible, e.g. when retroactively updating single values in the tables without updating the corresponding calculations.
Therefore, auditors often have to carry out a manual check of table correctness for each document.

During the manual auditing process,
one would apply knowledge about financial reports to evaluate which values correspond to which sums for each table
and recheck the correctness of the calculations.
This is of course highly time and labour intensive,
and mistakes are easy to make when auditing a large number of tables.

Automating this process however proves difficult.
While a machine has no problem evaluating the correctness of calculations for given table entries,
evaluating which values must sum up to which other values is a complex task.
Human auditors are either able to apply knowledge on financial reports or knowledge on structure of tables in general:
While tables are formatted in a way that human readers can easily evaluate how sums are structured, e.g. by sum values being at the bottom of the tables, indicted by text in the row headers, split from other values by bold lines,
teaching a machine to understand these intuitive rules is almost impossible.
A strict rule-based approach would be highly dependent on the formatting of specific tables and not generalize well. See Figure \ref{fig:sum_structure} for examples of easy and hard tables to parse.

A different approach to this problem is therefore ignoring the table structure altogether.
Treating single columns or the entire table as a (ordered) set of numbers,
one can try to find which values can be represented as a sum of a subset of all other numbers.
This problem can of course be solved exactly by deterministic algorithms.
However, the size of tables and magnitude of their entries (for financial documents) make many algorithmic approaches impractical.

In this paper,
we evaluate how stochastic algorithms based on gradient descent, minimizing the energy of Hopfield networks initialized with the problem formulation as parameters, can solve the problem of finding sums in large (both in size and magnitude) tables:
\begin{itemize}
    \item We restate the problem of finding sums in tables as the subset sum problem and briefly discuss known deterministic solving algorithms.
    \item We derive the general algorithm for gradient descent on Hopfield networks and how to restate the subset sum problem as a problem solvable by these networks (QUBOs).
    \item We evaluate our Hopfield algorithm on both artificial and real data and evaluate its performance in the context of financial auditing.
    \item We describe how quantum annealing hardware is able to solve QUBOs and therefore subset sum problems and run experiments to show that real quantum hardware is capable of solving subset sum problems.
\end{itemize}

% This paper is an extension of a preprint\footnote{\href{}{\texttt{(url removed for anonymous submission)}}}.
% We elaborate on the concepts presented in the preprint and additionally evaluate the application of quantum annealers and FPGA-hardware for the subset sum problem.

We provide code to reproduce all experiments in this paper at \url{https://github.com/fraunhofer-iais/quantum_subset_sum}.

\subsection{Subset Sum Formulation of the problem}

\begin{figure}[]
     \centering
     \begin{subfigure}[b]{0.45\columnwidth}
         \centering
         % \hspace*{1cm}       
        \includegraphics[height=3cm]{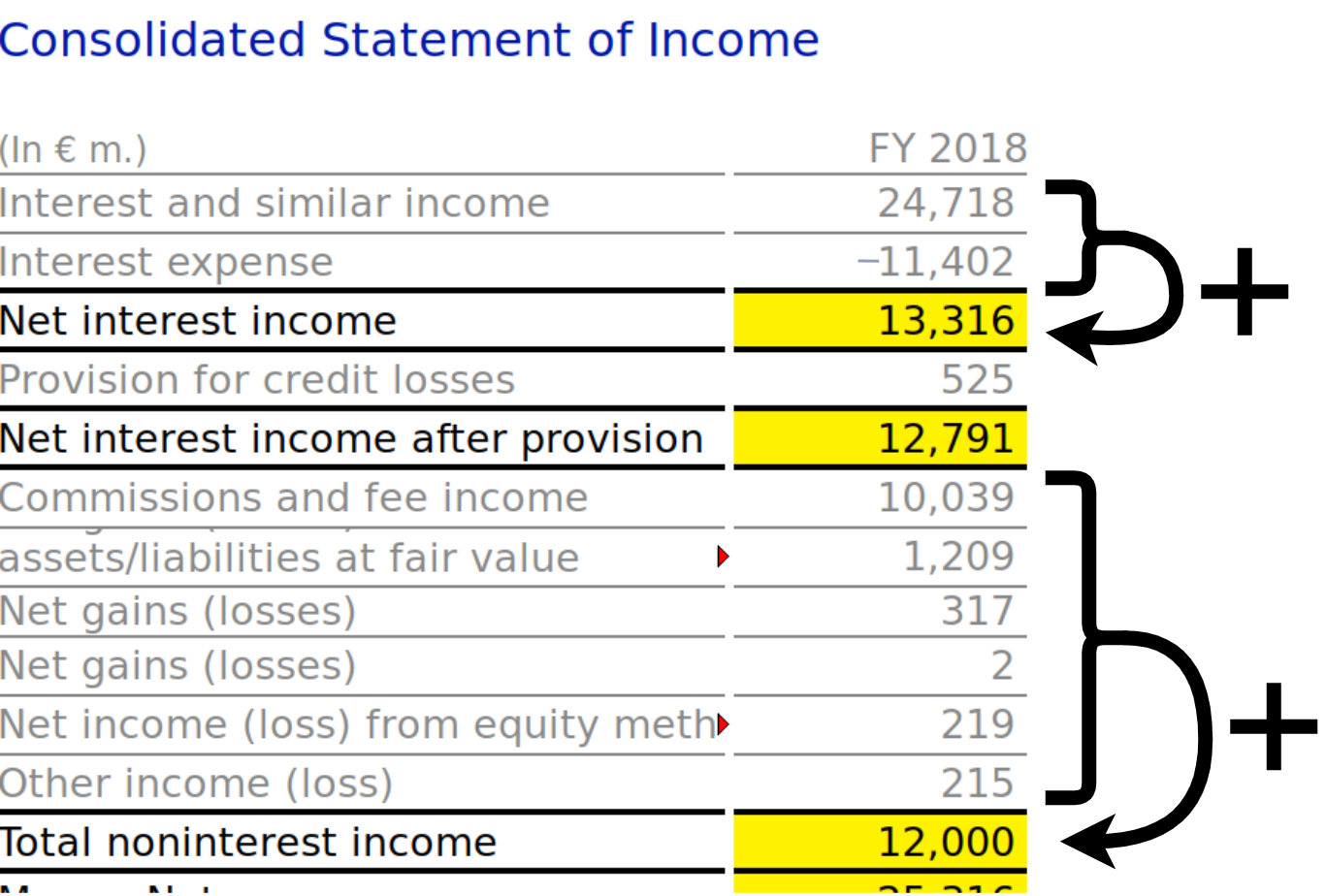}
         \label{fig:example_easy}
     \end{subfigure}%
     \hspace{0.5cm}
     \begin{subfigure}[b]{0.45\columnwidth}
         \centering
        \includegraphics[height=3cm]{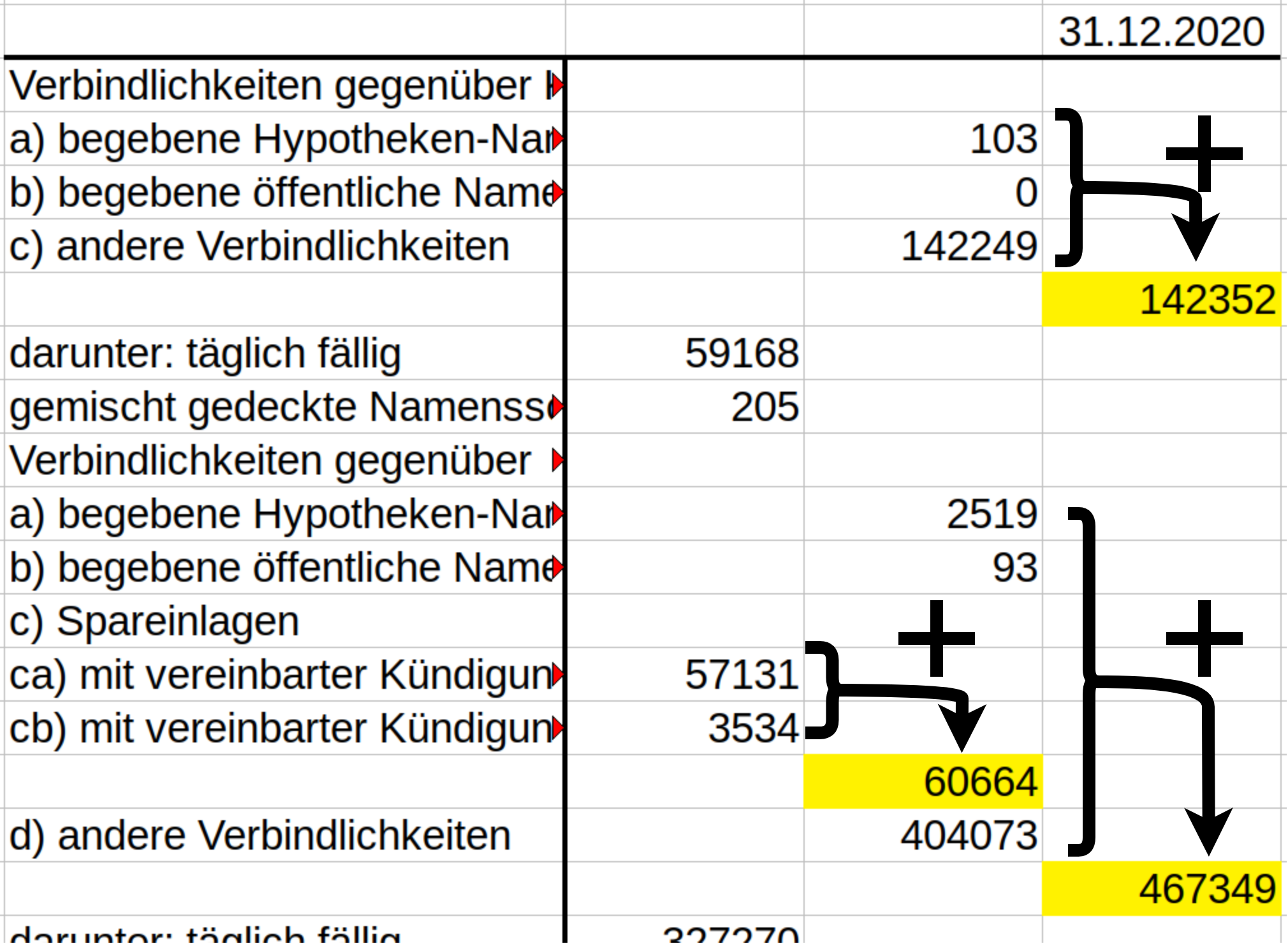}
         \label{fig:example_complicated}
     \end{subfigure}%
     \caption{Left: Example for easy sum structure. Column of a table containing financial performance indicators with annotated sums. Numbers from \cite{DBFDS}. Right: Example for more complicated sum structure. Sums extend over multiple columns. Numbers from \cite{DBKonzernabschluss}.}
    \label{fig:sum_structure}
\end{figure}

We call a set of rules that describe the behaviour of sums in a document, e.g. rows 1 and 2 sum up to row 3, rows 6 to 11 sum up to row 12, as a \emph{sum structure}.
See Figure \ref{fig:sum_structure} for two examples.

Many tables found in financial reports show the same sum structure in multiple columns, e.g. when comparing financial statements for several quarters and years.
Having an efficient algorithm for discovering sums in columns can aid consistency checks by applying the algorithm on one column,
extracting a sum structure and checking if the other columns also comply to the found sum structure.
If the new column does not comply to the same sum structure,
it is an indication for some inconsistency happening in the table.
This inconsistency can be displayed to an auditor for further investigation.

If no such repeated structure is present in the table (for example right side of Figure \ref{fig:sum_structure}),
an automatic parsing of sum structures still is valuable.
An application can extract all possible sums in the table and present its results to the auditor, 
who checks whether the parsed sum structure is reasonable, expected sums are missing or the displayed sums do not corespond to actual relations of table items.

Finding sum structures in tables is closely related to a well known problem in algorithmic combinatorics, the subset sum problem.
The subset sum problem is defined by a set of numbers
$\set{X} = \bigl\{ x_1, x_2, \ldots, x_n \bigr\} \subset \ZZ$ and a target sum $T \in \ZZ$.
We aim to find a subset $\set{Y}^* \subseteq \set{X}$,
such that the sum of the subset is equal to the target sum:
\begin{equation}
    \sum_{x \in \set{Y}^*} x = T.
\end{equation} 

In general, due to there being $2^n$ possible combinations of numbers for the subset, 
the problem is NP-hard.

In the framework of consistency checks and finding sums in tables,
we can consider the entire table as a set of numbers and apply the problem to each entry:
taken the entry as a target sum, is it possible to find a subset of all other numbers that sums up to the target?
Iterating a solving algorithm over each entry in the table yields a sum structure for the table.

While, there are several known algorithms for solving the subset sum problem,
The naive approach consists of cycling through each of the $2^n$ possible subsets for a total complexity of $\bigo(2^nn)$.
The algorithm can be improved by several heuristics and restrictions on the problem formulation (e.g. only positive integers) but the exponential complexity remains. \cite{exact_subset_sum}
Dynamic programming algorithms for solving the subset sum problem exactly in pseudo-polynomial time depend on the magnitude of values in $\set{X}$ \cite{pisinger_linear_time_algos}.
Large financial tables defy restrictions on both size and magnitude of values by containing many dozen individual entries and entries of up to billions of euros exact to one cent.

% \subsection{Classical solving algorithms and algorithms for approximate solutions}

% There are several known algorithms for solving the subset sum problem.
% The naive approach consists of cycling through each of the $2^n$ possible subsets, summing up all elements and comparing the sum to the target sum.
% This has a total complexity of $\bigo(2^nn)$.
% The algorithm can be improved by several heuristics and restrictions on the problem formulation (e.g. only positive integers) but the exponential complexity remains. \cite{exact_subset_sum}

% Additionally there exist dynamic programming algorithms for solving the subset sum problem exactly in pseudo-polynomial time, depending on the magnitude of values in $\set{X}$ \cite{pisinger_linear_time_algos}.

% Traditional algorithms can provide solutions in reasonable time for either a small amount of numbers or a small magnitude of values.
% Large financial tables defy both restriction by containing many dozen individual entries and entries of up to billions of euros exact to one cent.

\subsection{Rule-based algorithms for finding sums in financial tables}

The problem of finding sum structures in tables does not have to be broken down to the subset sum problem.
By ignoring the inherent structure and logic of the table, the complexity of the binary combination problem is increased.
Applying rule-based logic and understanding of the general structure of tables can result in efficient algorithms to solve the problem of finding sum structures in tables.

% These rules-based approaches can apply multiple heuristics to find sums:
% sums are generally more likely to be structured top-to-bottom,
% sums are likely to occur in the same row as the corresponding subset,
% the last entries in columns are likely to be sums.

However, rule-based approaches require specialization to each type of table.
See right side of Figure \ref{fig:sum_structure} for an example of a hard table to parse.
Numbers in one column can, but do not have to, correspond to sums in the next column.
Membership of a particular number to a column is determined by the accounting related properties of the corresponding value and any approach of parsing the internal table logic into strict rules is very unlikely to generalize to other tables.

\subsection{Related work}

Subset sum is a well known special case of the knapsack problem \cite{knapsack_problems}.
The run-time complexity of the subset sum problem depends on the parameters $n$, the number of input values, and $M$, the precision of the problem, i.e. the number of binary digits necessary to state the problem.
All deterministic exact solving algorithms are exponential in either $n$ or $M$. \cite{optimal_number_partitioning}
There exist dynamical programming approaches to subset sum \cite{faster_pseudopolynomial_time_algo}
which offer a polynomial run-time, 
however they require a space complexity of $\bigo(nT)$,
with $T$ being the target sum.
This space complexity is impractical for the proposed application in financial data analysis (see also Section \ref{sec:data}).

There is much work done on parallelization approaches and GPU-assisted computing to solve an adjusted problem statement in which one tries to maximize 
$\sum_{x \in \set{Y}} x$ under the restriction $\sum_{x \in \set{Y}} x \leq T$.
See \cite{low_space_algorithm_for_subset_sum} for a discussion of various algorithms.

In this work we reformulate the subset sum problem as a QUBO, quadratic unconstrained binary optimization problem.
QUBOs are closely related to Ising models and are a natural formulation for many related problems. \cite{ising_formulation}

We apply optimization over Hopfield networks \cite{hopfield_vector_quantization} to optimize an error term derived from the subset sum problem.
Hopfield networks are a special type of recurrent neural network which lends itself easily for solving QUBO problems. \cite{bauckhage_hopfield}
QUBO problems are also a main application of adiabatic quantum computing, where specialized quantum computers solve QUBOs by a process called quantum annealing. \cite{bauckhage_hopfield, dwave, quantum_annealing}.
In lieu of a quantum computer, specialized non-quantum hardware called field-programmable gate arrays (FPGA) are also able to very quickly and energy-efficiently solve QUBO problems \cite{Muecke2019-HAO}.

Lastly, our application lies in the field of automated financial analysis and auditing,
a field that has been growing quickly in recent years.
Specialized software, for example for text recommendation \cite{ali} and automatic anonymization \cite{anonymizer} are entering the market and we aim to add further to this toolbox.

\section{Subset-Sum as a QUBO}

A quadratic unconstrained binary optimization problem (QUBO) is defined by a function $f: \{0, 1\}^n \xrightarrow{} \RR$ which is a quadratic polynomial over its binary input variables,
\begin{equation}
    f(\vec{z}) = \sum_{ij} p_{ij} z_i z_j = \sum_{i \neq j} p_{ij} z_i z_j + \sum_i p_{ii} z_i.
\end{equation}
The QUBO problem consists of finding the optimal binary vector $\vec{z}^* \in \{0, 1\}^n$ such that
\begin{equation}
    \vec{z}^* = \argmin_{\vec{z}\in\{0, 1\}^n} f(\vec{z}).
\end{equation}

The problem can be rewritten in matrix notation as
  \begin{equation*}
  \vec{z}^* = \argmin_{\vec{z} \in \{0,1\}^n} \; \trn{\vec{z}} \mat{P} \, \vec{z} - \ipt{\vec{p}}{\vec{z}}
  \end{equation*}
  
  with a symmetric and hollow matrix $\mat{P} \in \RR^{n\times n}$
  and a vector $\vec{p} \in \RR^n$.
  
To convert the subset sum problem into a QUBO, we recall the problem statement.
Given set $\set{X} = \bigl\{ x_1, x_2, \ldots, x_n \bigr\}$ and target value $T$, determine a subset $\set{Y}^* \subseteq \set{X}$ such that $\sum_{x \in \set{Y}^*} x - T = 0$.
  
  The subset sum problem can therefore be stated as finding $\set{Y}^*$ such that
  \begin{equation}
  \set{Y}^*
  = \argmin_{\set{Y} \subseteq \set{X}} \left( T - \sum_{y \in \set{Y}} y \right)^2.
  \end{equation}

Collecting the numbers contained in set $\set{X}$ in a vector $\vec{x} = \trn{\bigl[ x_1, x_2, \ldots, x_n \bigr]} \in \mathbb{R}^n$
and introducing a binary indicator vector $\vec{z} \in \{0,1\}^n$ with entries
\begin{align*}
z_i & = 
  \begin{cases}
  1 & \text{if} \; x_i \in \set{Y} \\[1ex]
  0 & \text{otherwise}
  \end{cases}
\end{align*}
the subset sum problem can alternatively be written as
 \begin{align*}
  \vec{z}^*
  & = \amin{\vec{z} \in \{0,1\}^n} \; \bigl( \ipt{\vec{x}}{\vec{z}} - T \bigr)^2. 
  \end{align*}

 Expanding the equation we write
  \begin{align*}
  \vec{z}^*
  & = \argmin_{\vec{z} \in \{0,1\}^n} \; \ipt{\vec{z}}{\vec{x}}\,\ipt{\vec{x}}{\vec{z}} - 2\,T\, \ipt{\vec{x}}{\vec{z}} \textcolor{gray}{\,+\, const} \\[1ex]
  & \equiv \argmin_{\vec{z} \in \{0,1\}^n} \; \trn{\vec{z}} \mat{P} \, \vec{z} - \ipt{\vec{p}}{\vec{z}} \\
  &\eqqcolon \argmin_{\vec{z} \in \{0,1\}^n} E(\vec{z})
  \end{align*}
where we introduced the shorthands
\begin{equation} \label{def:Qq}
    \begin{split}
  \mat{P} & = \opt{\vec{x}}{\vec{x}}, \quad
  \vec{p} = 2\, T\, \vec{x}.
    \end{split}
\end{equation}

    Closely related to QUBOs are Ising Models, where we optimize over $\vec{s} \in \{-1, 1\}$ instead of $\vec{z} \in \{0, 1\}$:
    \begin{align*}
  \vec{s}^* &= \argmin_{\vec{s} \in \{-1,+1\}^n} \; \trn{\vec{s}} \mat{Q} \, \vec{s} + \ipt{\vec{q}}{\vec{s}}
  \end{align*}
  
  Both problem statements are in fact equivalent,
  with conversion via $\vec{z} = \tfrac{1}{2} (\vec{s} + \vec{1})$ and $\vec{s} = 2\vec{z} - \vec{1}$.
  
  Converting the QUBO derived from the subset sum problem above, we have
  \begin{align*}
  E ( \vec{z} )
  & = \trn{\vec{z}} \mat{P} \, \vec{z} - \ipt{\vec{p}}{\vec{z}} \\[1ex]
  & = \tfrac{1}{4} \, \trn{(\vec{s} + \vec{1})} \mat{P} \, (\vec{s} + \vec{1}) - \tfrac{1}{2} \, \ipt{\vec{p}}{(\vec{s} + \vec{1})} \\[1ex]
  & = \tfrac{1}{4} \, \trn{\vec{s}} \mat{P} \, \vec{s} + 
      \tfrac{1}{4} \, \trn{\vec{s}} \mat{P} \, \vec{1} + 
      \tfrac{1}{4} \, \trn{\vec{1}} \mat{P} \, \vec{s} - 
      \tfrac{1}{2} \, \ipt{\vec{p}}{\vec{s}} \textcolor{lightgray}{\,+\, const} \\[1ex]
  & = \tfrac{1}{4} \, \trn{\vec{s}} \mat{P} \, \vec{s} + 
      \tfrac{1}{2} \, \trn{\vec{1}} \mat{P} \, \vec{s} -
      \tfrac{1}{2} \, \ipt{\vec{p}}{\vec{s}} \textcolor{lightgray}{\,+\, const} \\[1ex]
  % & = \tfrac{1}{4} \, \trn{\vec{s}} \mat{P} \, \vec{s} +
  %     \tfrac{1}{2} \, \trn{\bigl( \mat{P} \, \vec{1} - \vec{p} \bigr)} \vec{s} \\[1ex]
  & \equiv \trn{\vec{s}} \mat{Q} \, \vec{s} + \ipt{\vec{q}}{\vec{s}} = E(\vec{s}),
  \end{align*}

with the shorthands
    \begin{equation}\label{def:Pp}
    \begin{split}
  \mat{Q} & = \tfrac{1}{4} \mat{P}, \quad
  \vec{q} = \tfrac{1}{2} \, \bigl( \mat{P} \, \vec{1} - \vec{p} \bigr),\end{split}
  \end{equation}
which holds since matrix $\mat{P} = \opt{\vec{x}}{\vec{x}} = \trn{\bigl( \opt{\vec{x}}{\vec{x}} \bigr)} = \trn{\mat{P}}$ is symmetric and we have $\trn{\vec{s}} \mat{P} \, \vec{1} = \trn{\vec{1}} \mat{P} \, \vec{s}$.

All in all, we can thus consider the subset sum problem as a minimzation problem over $\{-1, 1\}^n$ by
  \begin{equation} \label{def:min_bipolar}
  \vec{s}^* = \argmin_{\vec{s} \in \{-1,+1\}^n} \; \trn{\vec{s}} \mat{Q} \, \vec{s} + \ipt{\vec{q}}{\vec{s}}
  \end{equation}
  with $Q$ and $q$ defined in \eqref{def:Qq} and $P$ and $p$ defined in \eqref{def:Pp}.

\section{QUBO-Solving with Hopfield Networks}

A Hopfield Network is a recurrent neural net of
$n$ interconnected neurons.
The state of the network is described by a bipolar vector $\vec{s} \in \{-1, 1\}^n$.
Each neuron is connected to every other neuron, with connection weights given by a matrix $\mat{W} \in \RR^{n \times n}$. Moreover, each neuron $s_i$ is a bipolar threshold unit with threshold $\theta_i$, such that
\begin{equation}\label{eq:hopfield_update}
    s_i = \sign \left( \trn{\vec{w}_t} \vec{s} - \theta_i \right).
\end{equation}

A Hopfield network architecture is therefore fully described by a matrix $\mat{W}$, a vector $\vec{\theta}$ and a current state $\vec{s}$.
An update of the network is done via \eqref{eq:hopfield_update},
either for all neurons at once or only a subset of neurons.

Defining the energy of the Hopfield network in state $\vec{s}$ by
\begin{equation}
    E(\vec{s}) \coloneqq - \frac{1}{2} \trn{\vec{s}}\mat{W}\vec{s} + \trn{\vec{\theta}}\vec{s},
\end{equation}
we find that, if the weight matrix $\mat{W}$ is symmetric and is hollow (i.e. has diagonal of all zeros),
then the Hopfield energy can never increase when updating one neuron by \eqref{eq:hopfield_update}.
Since
\begin{equation}
\nabla E(\vec{s}) = -\mat{W}\vec{s} + \vec{\theta}
\end{equation}
the updates in \eqref{eq:hopfield_update} amount to $s_i = \sign(-\nabla E(\vec{s}))$ and each update performs gradient descent on $E(\vec{s})$.
As there are only $2^n$ possible states the network can be in, successive updates of single neurons will reach a local or global minimum after finitely many updates.
See \cite{bauckhage_hopfield} for details of the derivation.

This behaviour can be leveraged to solve problems stated as QUBOs.
Encoding the problem in weight and bias parameters $\mat{W}$ and $\vec{\theta}$,
such that minimum energy states
\begin{equation}
    \vec{s}^* = \argmin_{\vec{s} \in \{-1,+1\}^n} - \frac{1}{2} \trn{\vec{s}}\mat{W}\vec{s} + \trn{\vec{\theta}}\vec{s}
\end{equation}
solve the underlying QUBO problem,
the network may find solutions to the QUBO by the described gradient descent updates of single neurons.

To apply Hopfield networks to the subset sum problem, we
recall the problem statement as a minimization problem
over $\vec{s} \in \{-1, 1\}^n$ in \eqref{def:min_bipolar}.
Defining
\begin{equation}
    \begin{split}
        \mat{W} = -2 \mat{Q} = -\frac{1}{2} \mat{P}, \quad
        \vec{\theta} = \vec{q} = \frac{1}{2} (\mat{P}\vec{1} - \vec{p})
    \end{split}
\end{equation}
we find suitable weights and biases such that the state of a Hopfield network optimized to a global minimum encodes a solution to the subset sum problem.

Note that a Hopfield network with a random initialization does not necessarily converge to a global optimum and local optima are not solutions to the subset sum problem.
However, initializing the network multiple times with random states and running until convergence increases the chances of finding a global optimum.

The solving algorithm is easily extensible for efficient utilization of GPUs. See Algorithm \ref{algo:hopfield_gpu} for the GPU accelerated full solving algorithm.

\begin{algorithm}
\scriptsize
\caption{GPU-accelerated algorithm for solving the subset sum problem with Hopfield networks.}
\label{algo:hopfield_gpu}
\begin{algorithmic}
\State Given a vector of numbers $\vec{x} \in \ZZ^n$ and target sum $T \in \ZZ$,
\State construct Hopfield network weights $\mat{W}$ and biases $\vec{\theta}$ by
\State $\mat{W} \gets -\frac{1}{2} \vec{x}\trn{\vec{x}}, \quad \vec{\theta} \gets \frac{1}{2} \left( \vec{x}\trn{\vec{x}} \vec{1} - 2T\vec{x} \right)$.
\State Define batch size $m$ dependent on GPU memory.
\State Define maximum number of steps $t_{\text{max}}$, $t \gets 0$.
\While{No solution found and $t <  t_{\text{max}}$}
    \State Initialize batch of network states,
    \State sample $\mat{S} \in \{-1, 1\}^{\{m\times n\}}$
    \While{Hopfield energy has not converged}
        \State Calculate gradient: 
        \State $\nabla \vec{E} \gets - \mat{W}\vec{S} + \vec{1}\vec{\theta}$
        \State Select indices of maximum change:
        \State $\vec{I}_j \gets \argmin_{i} \nabla \vec{E}_{ji}$ 
        \State Update $\vec{S}_j$ at index $\vec{I}_j$: 
        \State $\vec{S}_{j\vec{I}_j} \gets \sign(\trn{\mat{W}_{\vec{I}_j}}\vec{S}_j - \vec{\theta}_{\vec{I}_j})$
        \State Calculate Hopfield energies of each state in the batch:
        \State $\vec{E}_j \gets -\mat{S}\mat{W}\mat{S}_{jj} + (\trn{\vec{\theta}}\vec{S})_j$
    \EndWhile
    \For{$j \in \{1, \dots, m\}$}
    \If{Hopfield energy is minimal: $\vec{E}_j = 0$}
        \State End algorithm. 
        \State \Return $\{\vec{x}_i | \vec{S}_{ji} = 1 \}$
    \EndIf
    \EndFor
    \State $t \gets t + 1$
    \EndWhile
\State Algorithm has not found a solution in $s_{\text{max}}$ steps.
\State End algorithm.
\State \Return None.
\end{algorithmic}
\end{algorithm}

% \begin{figure}
%     \centering
%     \includegraphics[width=0.8\columnwidth]{figures/fpga.jpg}
%     \caption{Picture of the special purpose hardware }
%     \label{fig:fpga}
% \end{figure}

\section{Experiments on Hopfield Networks}

\subsection{Data}
\label{sec:data}

We conduct experiments with both artificial data and real data.
To create artificial data we uniformly sample $n$ integers between $X_{\text{min}}$ and $X_{\text{max}}$,
select $k$ of the sampled integers at random and calculate the target sum $T$ as the sum of the selected integers.

% Note that the selection of $n$, $X_{\text{min}}$ and $X_{\text{max}}$ defines the number of possible solutions to the problem and therefore influences the difficulty of finding a solution.
% Given a set $\set{X}$ of $n$ integers between $X_{\text{min}} < 0$ and $X_{\text{max}} > 0$, the sum of any subset of $\set{X}$
% must be in the interval 
% \begin{equation}
%     \set{T} = \left[nX_{\text{min}}, nX_{\text{max}}\right]
% \end{equation}
% for a total of $\# \set{T} = n\left(X_{\text{max}} - X_{\text{min}}\right)$.
% However, there are $2^n$ possible subsets of $\set{X}$.
% For many combinations of $n$, $X_{\text{min}}$ and $X_{\text{max}}$ we have
% $2^n \gg n\left(X_{\text{max}} - X_{\text{min}}\right)$ and therefore some target values must have multiple solutions.
% Finding solutions for a problem with many distinct solutions is of course easier than finding one correct solution in $2^n$ possible combinations.

We construct artificial data for $n = 16, 32, 64, 128, 256$ ($k = 4, 8, 8, 8, 8$ respectively) and values in $\left[-X_{\text{max}}, X_{\text{max}}\right]$ for $X_{\text{max}} = 1\mathrm{e}{+4}, 1\mathrm{e}{+5}, 1\mathrm{e}{+6}$.
% See Table \ref{tab:art_data_configs} for an overview of $R = 2^n / 2nX_{\text{max}}$,
% the expected number solutions when sampling a set of numbers and a random target solution.
For each configuration, we sample $M = 5$ different subset sum problems.

We evaluate our algorithm on a set of real data problems.
We parse a financial report \cite{DBFDS}
containing multiple sheets with financial reports for the quarters from Q1 2019 to Q4 2020 for a total of 190 individual subset sum problems.
Each column contains numbers describing amounts up to multiple billion euro, exact to one cent, for a total of 14 significant figures. See Table \ref{tab:fin_data_configs} for details on the dataset.
In contrast to the artificial dataset, the real dataset only contains exactly one solution for most subset sum problems.

\begin{table}[]
    \centering
    \scriptsize
    \begin{tabular}{l|rrrrrr}
    \toprule
    name & $M$ & n & $X_{\text{min}}$ & $X_{\text{max}}$  \\ %& $R$\\
    \midrule
    assets       & 15 & 17 & $1.7\mathrm{e}{+07}$ & $1.5\mathrm{e}{+14}$  \\ %& $5.1\mathrm{e}{-11}$
    consincome   & 49 & 31 & $-2.5\mathrm{e}{+12}$ & $2.5\mathrm{e}{+12}$ \\ %& $1.4\mathrm{e}{-05}$
    liabilities  & 29 & 20 & $2.5\mathrm{e}{+10}$ & $1.5\mathrm{e}{+14}$  \\ %& $3.5\mathrm{e}{-10}$
    net revenues & 97 & 25 & $-5.3\mathrm{e}{+10}$ & $2.5\mathrm{e}{+12}$ \\ %& $5.3\mathrm{e}{-07}$
    \bottomrule
    \end{tabular}
    \caption{Configurations of parsed financial data. Each group describes one table in the financial report, 
    with $M$ individual subset sum problems (i.e. sums described in the table), a column length of $n + 1$ which corresponds to $n$ values in the subset sum problem, and values between $X_{\text{min}}$ and $X_{\text{max}}$.
    % $R$ describes the expected number solutions when sampling a set of numbers and a random target solution defined by $n$, $X_{\text{min}}$ and $X_{\text{max}}$. We see that each $R \ll 1$ and in most columns there is only one unique solution to the subset sum problem.
    }
    \label{tab:fin_data_configs}
\end{table}

\subsection{Experiments and Results}

\begin{figure}[]
     \centering
     \begin{subfigure}[b]{0.4\textwidth}
         \centering
         \hspace{-.8cm}\includegraphics[height=4cm]{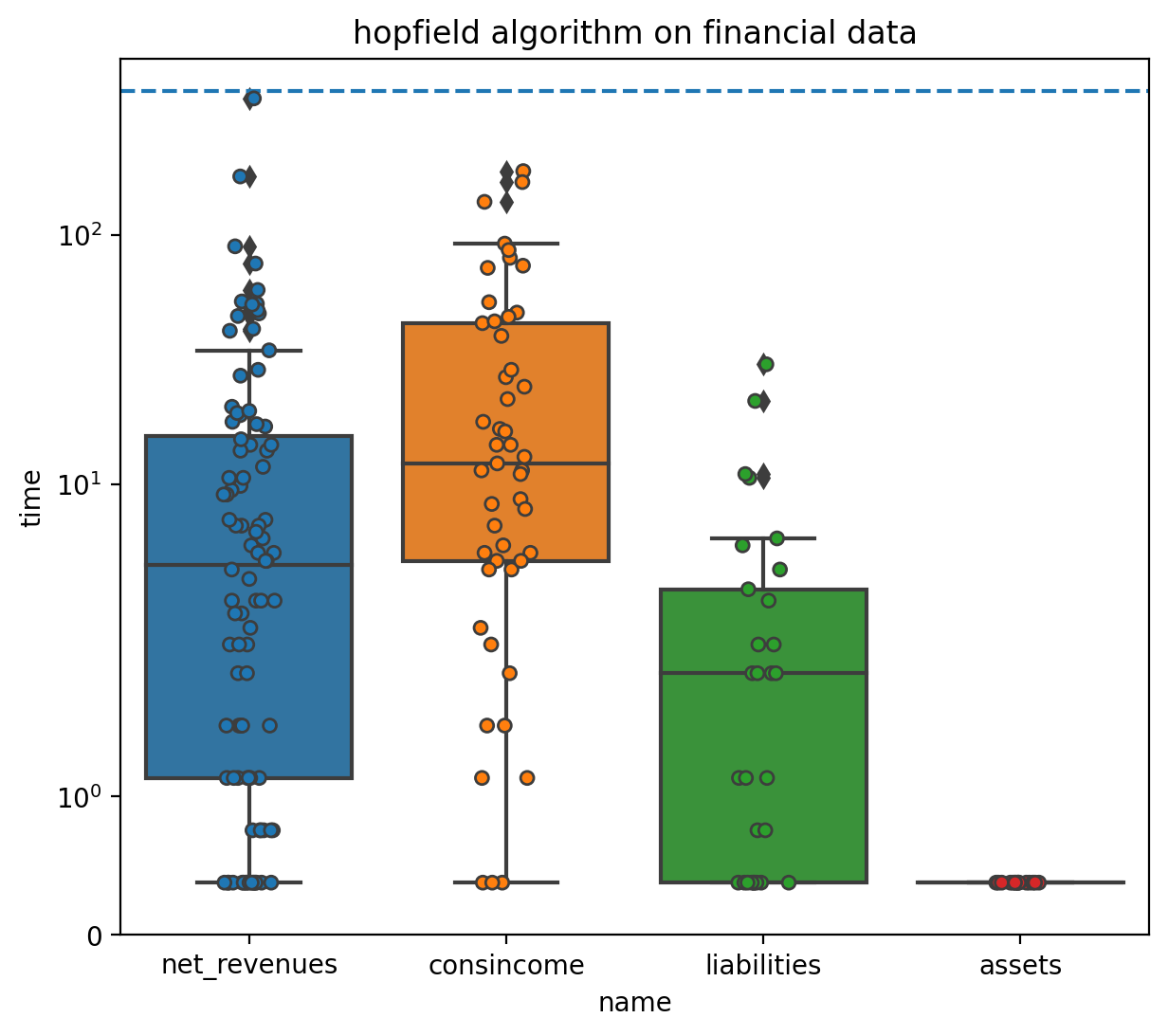}
         \label{fig:plot_count_artificial}
     \end{subfigure}%
     
     \begin{subfigure}[b]{0.4\textwidth}
         \centering
          \hspace{3cm}\includegraphics[height=4cm]{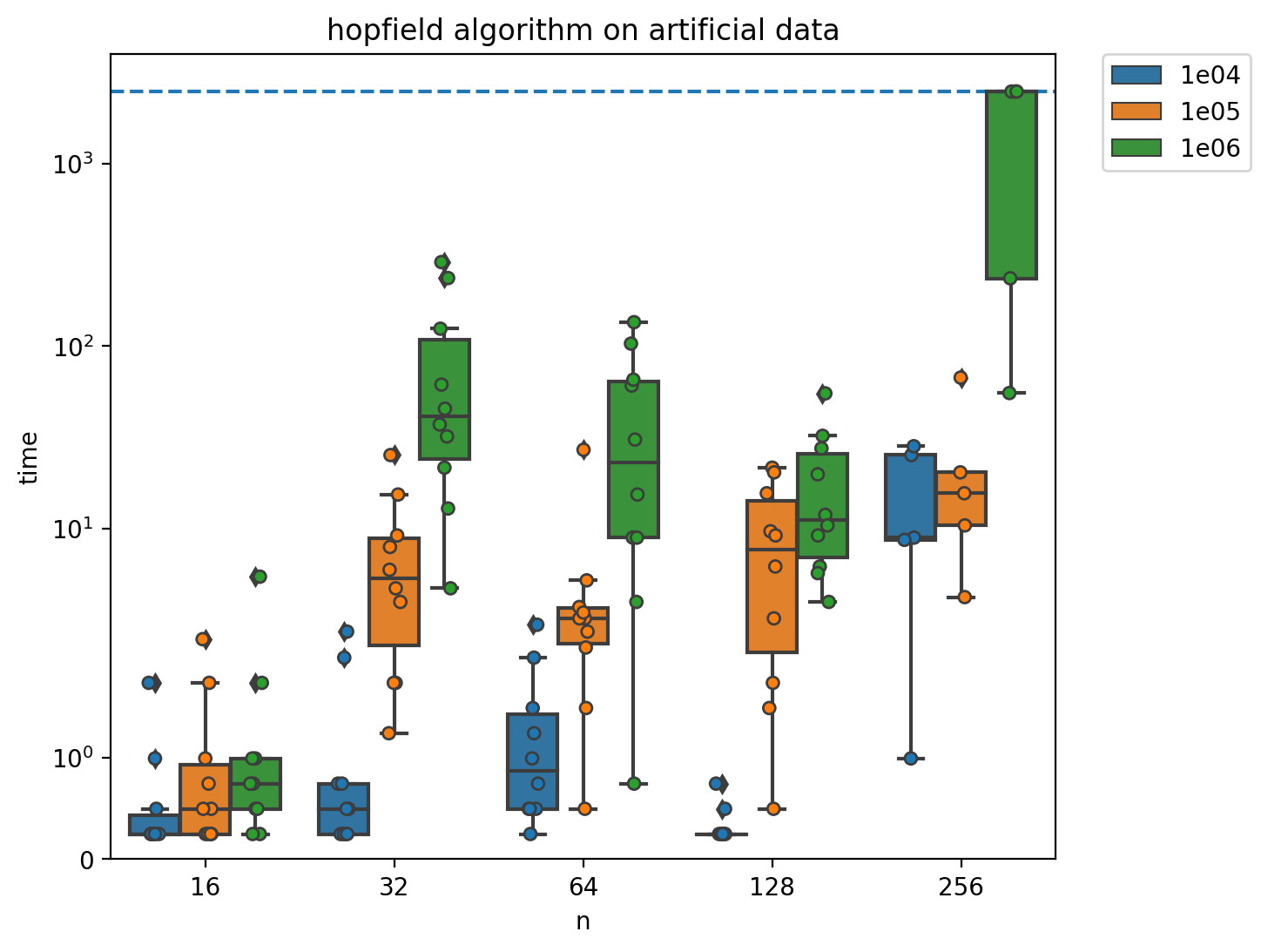}
         \label{fig:plot_hopfield_artificial}
     \end{subfigure}%
    %  \hfill
     \caption{Benchmarking the Hopfield algorithm on the parsed financial documents (top) and artificial data (bottom). The dotted line roughly indicates the cutoff for computation.}
     \label{fig:hopfield_plots}
\end{figure}

We run the Hopfield algorithm on artificial data for up to $1\mathrm{e}{+8}$ initializations in batches of $1\mathrm{e}{+4}$ on one NVIDIA A100 GPU, for a maximum computation time of around 8 minutes. 
For almost all configurations the algorithm reliably finds a correct solution for all samples.
Only for the configuration of $n = 256$ and $X_{\text{max}} = 1\mathrm{e}{+6}$ not all samples are solved in the specified maximum number of runs (2 of 5 found).
See Figure \ref{fig:hopfield_plots} (bottom) for comparison of computation time against $n$ and  $X_{\text{max}}$. We see that the number of values $n$ has a smaller impact on the computation time until a solution is found than the magnitude of the numbers  $X_{\text{max}}$.

We run the Hopfield algorithm with the same configuration of the financial data.
The algorithm finds a correct solution to the problem under the maximum number of iterations in all cases.
Note that unlike for most of the artificial data, the combination of $n$ and $X_{\text{max}}$ lead to a situation where each problem likely only contains one correct solution, which is found by our algorithm. 
See Figure \ref{fig:hopfield_plots} (left) for a comparison of computation time for different tables.
We see that average computation time clearly increases with the size of the problem, i.e. amount of values in the table.

See Tables \ref{tab:results_hopfield_artificial} and \ref{tab:results_hopfield_financial} for additional statistics on the runs.
In total, we conclude that optimization of binary vectors with Hopfield networks is a reliable algorithm for solving subset sum problems and can be applied to real-world examples.

\begin{table}[]
\centering
\scriptsize
    \begin{tabular}{ll|rr}
    \toprule
    n & $X_{\text{max}}$ &  $t_\text{mean}$           &     $n_\text{mean}$          \\
    \midrule
    16  & $1\mathrm{e}{+4}$ &            0.4 &  $2.0\mathrm{e}{+4}$ \\
        & $1\mathrm{e}{+5}$ &            0.7 &  $3.2\mathrm{e}{+4}$ \\
        & $1\mathrm{e}{+6}$ &            1.1 &  $4.9\mathrm{e}{+4}$ \\
    32  & $1\mathrm{e}{+4}$ &            0.7 &  $3.3\mathrm{e}{+4}$ \\
        & $1\mathrm{e}{+5}$ &            6.9 &  $3.1\mathrm{e}{+5}$ \\
        & $1\mathrm{e}{+6}$ &           77.2 &  $3.4\mathrm{e}{+6}$ \\
    64  & $1\mathrm{e}{+4}$ &            1.0 &  $4.5e\mathrm{+}{4}$\\
        & $1\mathrm{e}{+5}$ &            4.9 &  $2.1\mathrm{e}{+5}$ \\
        & $1\mathrm{e}{+6}$ &           40.1 &  $1.7\mathrm{e}{+6}$ \\
    % 128 & $1\mathrm{e}{+04}$ &            0.3 &  $1.3e\mathrm{+}{04}$\\
    %     & $1\mathrm{e}{+05}$ &            9.0 &  $3.6\mathrm{e}{+05}$ \\
    %     & $1\mathrm{e}{+06}$ &           18.4 &  $7.3\mathrm{e}{+05}$ \\
    % 256 & $1\mathrm{e}{+04}$ &            14.4 &  $1.3e\mathrm{+}{04 }$\\
    %     & $1\mathrm{e}{+05}$ &            23.5 &  $3.6\mathrm{e}{+05}$ \\
    %     & $1\mathrm{e}{+06}$ &           145.44 &  $7.3\mathrm{e}{+05}$ \\
    \bottomrule
    \end{tabular}
    \begin{tabular}{ll|rr}
    \toprule
    n & $X_{\text{max}}$ &  $t_\text{mean}$           &     $n_\text{mean}$          \\
    \midrule
    % 16  & $1\mathrm{e}{+04}$ &            0.4 &  $2.0\mathrm{e}{+04}$ \\
    %     & $1\mathrm{e}{+05}$ &            0.7 &  $3.2\mathrm{e}{+04}$ \\
    %     & $1\mathrm{e}{+06}$ &            1.1 &  $4.9\mathrm{e}{+04}$ \\
    % 32  & $1\mathrm{e}{+04}$ &            0.7 &  $3.3\mathrm{e}{+04}$ \\
    %     & $1\mathrm{e}{+05}$ &            6.9 &  $3.1\mathrm{e}{+05}$ \\
    %     & $1\mathrm{e}{+06}$ &           77.2 &  $3.4\mathrm{e}{+06}$ \\
    % 64  & $1\mathrm{e}{+04}$ &            1.0 &  $4.5e\mathrm{+}{04}$\\
    %     & $1\mathrm{e}{+05}$ &            4.9 &  $2.1\mathrm{e}{+05}$ \\
    %     & $1\mathrm{e}{+06}$ &           40.1 &  $1.7\mathrm{e}{+06}$ \\
    128 & $1\mathrm{e}{+4}$ &            0.3 &  $1.3e\mathrm{+}{4}$\\
        & $1\mathrm{e}{+5}$ &            9.0 &  $3.6\mathrm{e}{+5}$ \\
        & $1\mathrm{e}{+6}$ &           18.4 &  $7.3\mathrm{e}{+5}$ \\
    256 & $1\mathrm{e}{+4}$ &            14.4 &  $1.3e\mathrm{e}{+4}$\\
        & $1\mathrm{e}{+5}$ &            23.5 &  $3.6\mathrm{e}{+5}$ \\
        & $1\mathrm{e}{+6}$ &           145.44 &  $7.3\mathrm{e}{+5}$ \\
    \bottomrule
    \\
    \\
    \\
    \end{tabular}
    \caption{Results for the Hopfield algorithm on artificial data configurations. 
    The algorithm found a solution for every given sample.
    For each data configuration we describe the mean time to find a solution in seconds and the mean number of individual runs until a solution is found.}
    \label{tab:results_hopfield_artificial}
\end{table}

\begin{table}[]
    \centering
    \scriptsize
    \begin{tabular}{l|rrrrrr}
    \toprule
    name             & $n$ & $M$ & $M_{\text{found}}$ &  n  & mean time &         mean runs    \\
    \midrule
        assets       & 17 & 15 & 15 & 17  & 0.4&  $1.0\mathrm{e}{+04}$ \\
        consincome   & 31 & 49 & 49  & 31 & 30.4 &  $8.0\mathrm{e}{+05}$ \\
        liabilities  & 20 & 29 & 29 & 20  & 4.1 &  $1.1\mathrm{e}{+05}$ \\
        net revenues & 25 & 97 & 97 & 25 & 17.6 &  $4.6\mathrm{e}{+05}$ \\
    \bottomrule
    \end{tabular}
    \caption{Results for the Hopfield algorithm on artificial data configurations. 
    For each data configuration we decribe the mean time to find a solution in seconds and the mean number of individual runs until a solution is found.
    Of $M$ subset sum problems on each table we solved $M_{\text{found}}.$}
    \label{tab:results_hopfield_financial}
\end{table}

\section{QUBO-Solving with Quantum Annealers}

The runtime of all algorithms presented in this paper so far still depends on the complexity of the problem in terms of the amount of numbers in the table.
This limitation can be partially erased by application of quantum computing hardware.

A quantum annealer is a type of adiabatic quantum computer that applies the physical quatum annealing process to solve problems.
A quantum annealer consists of multiple qubits, which like regular bits can be in the binary states 0 and 1, but also in a superposition between the states.

Solving QUBOs comes naturally to a quantum annealing machine \cite{bauckhage_hopfield}.
Each qubit is initialized with an magnetic field called a bias, and multiple qubits can be linked together such that they influence each other, which is called coupling.
To solve QUBOs, each qubit represents one binary variable $z_i$, the bias is initialized by vector $\vec{p}$ in \eqref{def:Pp} and the coupling weights are given by matrix $\mat{P}$ in \eqref{def:Pp}.
The initialized system of quantum objects describes an energy landscape and the qubits will naturally find a state of low or minimal energy.
At the end of the quantum annealing process, each qubit will collapse into a binary state of either 0 or 1.
If the system found a global minimum of the energy landscape,
the resulting set of collapsed qbits describes the solution to the QUBO.
Due to the physical nature of this process, in theory quantum annealers can find the optimal solution to problems of arbitrary size almost instantly.

However, due to the structure of the graph connecting the qubits,
a total of $\bigo(M)$ physical qubits can only describe $\bigo(\sqrt{M})$ logical qubits.
The most prominent manufacturer of quantum annealers is D-Wave Systems, 
which provides machines with a maximum of \num{5436} physical qubits.
Meaning the largest quantum annealer available today would be able to solve QUBOs of $70$ variables.
A smaller machine with $\num{2000}$ qubits could solve problems of size $50$ \cite{dwave_bits}.

We test the capability of current quantum annealers on a small test dataset. In a table \cite{adidas_table_small} with 17 entries of values between $X_\text{min} = -2.8\mathrm{e}{+04}$ and $X_\text{max} = 4.3\mathrm{e}{+06}$ we identify 6 individual subset sum problems.
We convert each problem into a QUBO and run the proposed algorithm on a DWave Advantage quantum annealer with 5436 physical qubits.
Table \ref{tab:quantum_annealer_results} describes the results of this experiment.
We see that the quantum computer finds the correct proposed solution for each problem.

% Additionally, the algorithm often finds multiple solutions that solve the underlying subset sum problem,
% since most sum structures have multiple solutions:
% a number $x$ that belongs to a certain sum can itself be a sum of other numbers $x_1, \dots, x_n$. In this case, both the solutions containing $x$ and containing $x_1, \dots, x_n$ are correct. See for example Figure \ref{fig:sum_structure} (left), in which rows $\{1, 2\}$ sum up to 3, and therefore both $\{1, 2, 4\}$ and $\{3, 4\}$ sum up to 5.
% Any software applying subset sum solving algorithms must keep this in mind when presenting solutions to auditors in a clearly structured way.

We run an additional experiment with subset sum problems from a table \cite{adidas_table_large} with 47 entries of $X_\text{min} = -9.0\mathrm{e}{+02}$ and $X_\text{max} = 2.6\mathrm{e}{+03}$.
While the quantum architecture is able to correctly process each QUBO and find solutions for each subset sum problem,
the larger amount of numbers with small magnitude leads to issues described in section \ref{sec:data}, where the number of individual numbers leads to multiple solutions of each target sum, most of which do not correspond to any real-world relationships.

For the time being, quantum annealers are very costly to operate and not entirely practical for this application.
However, as technology advances, quantum computing becomes more and more affordable and practical which 
Once quantum annealing is ready for application on the scale of problems described in this paper,
one can directly apply the methodology described in this paper to apply quantum hardware in the auditing process.

\begin{table}[]
    \centering
    \begin{tabular}{c|rrrr}
    \toprule
        Problem & $N_{\text{o}}$     & $N_{\text{c}}$ & found \\
    \midrule
        \#1     & 1                      & 1                    & yes  \\  
        \#2     & 5                      & 1                    & yes  \\
        \#3     & 3                      & 2                    & yes  \\
        % \#4     & 2                      & 2                    & yes  \\
        % \#5     & 8                      & 6                    & yes  \\
        % \#6     & 5                      & 5                    & yes  \\   
    \bottomrule
    \end{tabular}
    \hspace{0.5cm}
    \begin{tabular}{c|rrrr}
    \toprule
        Problem & $N_{\text{o}}$     & $N_{\text{c}}$ & found \\
    \midrule
        % \#1     & 1                      & 1                    & yes  \\  
        % \#2     & 5                      & 1                    & yes  \\
        % \#3     & 3                      & 2                    & yes  \\
        \#4     & 2                      & 2                    & yes  \\
        \#5     & 8                      & 6                    & yes  \\
        \#6     & 5                      & 5                    & yes  \\   
    \bottomrule
    \end{tabular}
    \caption{Results of the experiments run on a DWave Advantage quantum annealing computer.
    Each problem stems from a parsed table \cite{adidas_table_small} with $N = 17$ entries of smaller magnitude between $X_\text{min} = -2.8\mathrm{e}{+04}$ and $X_\text{max} = 4.3\mathrm{e}{+06}$.
    We evaluate the number of global optima of the energy landscape the annealing process finds ($N_{\text{optim}}$),
    the number of optima that correspond to a correct solution of the subset sum problem ($N_{\text{correct}}$),
    and whether one of the proposed solutions corresponds to the correct target solution.
    }
    \label{tab:quantum_annealer_results}
\end{table}

\section{Conclusion and Outlook}

In this work we investigated how the subset sum problem plays a vital part in the automation of the financial auditing process and
how the subset sum problem can be restated as a well known problem architecture
which can be solved by the application of gradient descent on the energy landscape of Hopfield networks.
We found that the proposed algorithm reliably finds correct sum structures for artificial and real data.

We gave an overview of the capabilities of adiabatic quantum computers for the proposed task and its current limitations.
We evaluated the capability of quantum annealers for the subset sum problem and found that for problems with small range of values the algorithm reliable finds correct solutions.
As quantum hardware becomes more capable, it will make the use of quantum computing for real applications in finance a realistic possibility.  
% We found that some of the discussed limitations of quantum computers can be overcome by the application of special purpose hardware for QUBO-solving. 
% Future work will include evaluation of FPGAs architectures for solving subset sum problems on the given datasets.

In the near future, the algorithm will be ready to deploy on existing smart auditing software to directly benefit auditors in their daily work.

\section{Acknowledgements}
In parts, the authors of this work were supported by the Smart Data Innovation Challenges (SDI-C) Project ``Solving Accounting Optimization Problems in the Cloud'', which was funded by the Federal Ministry of Education and Research (BMBF) of Germany.

In parts, this research has been funded by the Federal Ministry of Education and Research of Germany and the state of North-Rhine Westphalia as part of the Lamarr-Institute for Machine Learning and Artificial Intelligence, LAMARR22B.

\bibliographystyle{ieeetr}
\bibliography{bibliography}

\end{document}